\documentclass[11pt,leqno]{article}
\usepackage{a4}
\usepackage[T1]{fontenc}
\usepackage[english]{babel}
\usepackage{latexsym}
\usepackage{amsmath}
\usepackage{amssymb}
\usepackage{mathrsfs}
\usepackage{color}
\usepackage{cite}
\usepackage{titling}
\usepackage{bbold}
\makeatletter
\@addtoreset{equation}{section} 
\makeatother
\newcommand{\ud}{\mathrm{d}}
\newcommand{\de}{\mathrm{D}}
\newcommand{\p}{\mathrm{p}}

\newcommand{\re}{\mathrm{Re}}
\newcommand{\im}{\mathrm{Im}}
\newcommand{\e}{\mathrm{e}}
\newcommand{\id}{\mathrm{Id}}
\newcommand{\R}{\mathbb{R}}
\newcommand{\N}{\mathbb{N}}

\newcommand{\Q}{\mathbb{Q}}

\newcommand{\ds}{\displaystyle}
\newcommand{\vfi}{\varphi}

\newenvironment{pr}{\vspace{5pt}\textbf{{\small Proof :}}\\}{\hspace{\stretch{1}}\rule{1ex}{1ex}\vspace{5pt}}
\newtheorem{thm}{Theorem}[section]

\usepackage{fancyhdr}
\fancypagestyle{plain}{\fancyhead{}}
\title{Remark on the pointwise stabilization of an elastic string equation}
\author{{FATHI HASSINE}\\\\ \textit{UR Analysis and Control of PDE 13ES64}\\ \textit{Department of Mathematics, Faculty of Sciences of Monastir}\\ \textit{University of Monastir, 5019 Monastir, Tunisia}\\ \textit{email:} \texttt{fathi.hassine@fsm.rnu.tn}}
\date{}
\begin{document}
\maketitle
\begin{center}
\abstract{
We consider an initial and boundary value problem the one dimensional wave equation with damping concentrated at an interior point. We prove a result of a logarithmic decay of the energy of a system with homogeneous Dirichlet boundary conditions. The method used is based on the resolvent estimate approach which derives from the Carleman estimate technique. Under an algebraic assumption describing the right location of the actuator, we prove a logarithmic decay of the energy of solution. We show that this assumption is lower than the one given by~\cite{Tuc} and~\cite{AHT} which depends on the diophantine approximations properties of the actuator's location.
}
\end{center}


\vspace{10pt}
\section{Introduction}
In recent years a lot of papers were devoted to the study of elastic structures with pointwise stabilizer:
\begin{enumerate}
	\item Placed on the boundary where several authors~\cite{C},~\cite{L} and~\cite{QR} have shown that the energy has the best design that satisfies what's called the uniform exponential stabilization property, i.e  there exist $M\geq1$ and $\mu>0$ such that the energy of the system decay exponentially
	$$
	E(t)\leq M.E(0)\e^{-\mu t},\;\forall\,t>0.
	$$
	\item Placed inside the span where several authors~\cite{CCW}~,\cite{Ho},~\cite{JTZ} and~\cite{Tuc} have shown at first strong stabilization depending on the position of the actuator and secondly they have shown uniform and non-uniform exponential decay of energy depending on the boundary conditions.
\end{enumerate}
The main purpose of the present paper is to study the stabilization model of a one dimensional coupled wave equation system with pointwise damping that models the vibrations of a string. The structure is formed by two coupled vibrating strings. It can be conceived as segments of power transmission lines, aerial cable/railway systems or the upper cable part of an idealized suspension bridge. More precisely we consider the following initial boundary value problem
\begin{eqnarray}\label{p1}
\ddot{u}_{1}(x,t)-u_{1}''(x,t)=0&\text{in}&(0,\xi)\times(0,+\infty),
\\
\ddot{u}_{2}(x,t)-u_{2}''(x,t)=0&\text{in}&(\xi,1)\times(0,+\infty),
\end{eqnarray}
where $u_{1}$ and $u_{2}$ denote the transverse displacement at the point $x$ and time $t$, with the two coupling transmission conditions at the interior node $\xi\in(0,1)$,
\begin{eqnarray}\label{p2}
u_{1}(\xi,t)=u_{2}(\xi,t)&\text{for}&t\in(0+\infty),
\\
u_{1}'(\xi,t)-u_{2}'(\xi,t)=\dot{u}_{1}(\xi,t)&\text{for}&t\in(0+\infty),
\end{eqnarray}
that describes the continuity of displacement for the first and the discontinuity of vertical force component for the second, besides Dirichlet boundary conditions at the right and left ends $x=0$ and $x=1$ are considered here
\begin{eqnarray}
u_{1}(0,t)=u_{2}(1,t)=0&\text{for}&t\in(0+\infty).
\end{eqnarray}
where we recall that the prime is the space derivative and the dot is the time derivative, and that the initial boundary conditions are given by
\begin{eqnarray}
u_{1}(x,0)=u_{1}^{0}(x),\;\dot{u}_{1}(x,0)=u_{1}^{1}(x)&\text{for}&x\in(0,\xi),
\\
u_{2}(x,0)=u_{2}^{0}(x),\;\dot{u}_{2}(x,0)=u_{2}^{1}(x)&\text{for}&x\in(\xi,1).\label{p3}
\end{eqnarray}

If $u_{1}$ and $u_{2}$ are the solution of~\eqref{p1}-\eqref{p3} we define the energy for $(u_{1},u_{2})$ at instant $t$ by
\[
\mathcal{E}(u_{1},u_{2})(t)=\frac{1}{2}\left(\int_{0}^{\xi}|\dot{u}_{1}(x,t)|^{2}+|u_{1}'(x,t)|^{2}\,\ud x+\int_{\xi}^{1}|\dot{u}_{2}(x,t)|^{2}+|u_{2}'(x,t)|^{2}\,\ud x\right).
\]
A simple formal calculation shows that the energy of the solution of~\eqref{p1}-\eqref{p3} is decreasing over the time as follows
\[
\mathcal{E}(u_{1},u_{2})(t_{1})-\mathcal{E}(u_{1},u_{2})(t_{2})=\int_{t_{1}}^{t_{2}}|\dot{u}_{1}(\xi,t)|^{2}\,\ud t,\;\;\forall\,t_{1},\,t_{2}\geq 0.
\]
Noting by $X=H_{0}^{1}(0,1)\times L^{2}(0,1)$ the space embedding with the norm
\[
\|(u,v)\|^{2}=\|u'\|_{L^{2}(0,1)}^{2}+\|v\|_{L^{2}(0,1)}^{2}
\]
in which we define the operator $\mathcal{A}$ by
\[
\mathcal{A}=\left[\begin{array}{cc}
0&\id
\\
\ds\frac{\ud^{2}}{\ud x^{2}}&0
\end{array}\right]
\]
with domain
\begin{equation*}
\begin{split}
\mathcal{D}(\mathcal{A})=\{(u,v)\in H_{0}^{1}(0,1)\times H_{0}^{1}(0,1):u_{|(0,\xi)}\in H^{2}(0,\xi),\,u_{|(\xi,1)}\in H^{2}(\xi,1),
\\
u'(\xi^{+})-u'(\xi^{-})=v(\xi)\}.
\end{split}
\end{equation*}
By noting 
$$u=\left\{\begin{array}{ll}
u_{1}&\text{in }(0,\xi)
\\
u_{2}&\text{in }(\xi,1),
\end{array}\right.
$$
system~\eqref{p1}-\eqref{p3} can be rewritten as the following Cauchy problem
\[
\left(\begin{array}{c}
\dot{u}
\\
\dot{v}
\end{array}\right)=\mathcal{A}\left(\begin{array}{c}
u
\\
v
\end{array}\right),\qquad\forall(u,v)\in\mathcal{D}(\mathcal{A}).
\]
It's well known (see for instance~\cite{CCW}) that for every $\xi\in (0,1)$ the operator $\mathcal{A}$ generates a $C_{0}$-semigroup of contraction. Then the problem~\eqref{p1}-\eqref{p3} is well-posed where the solution $u$ satisfies
\[
u\in\mathrm{C}([0,+\infty),\mathcal{D}(\mathcal{A}))\cap\mathrm{C}^{1}([0,+\infty),X)
\]
if the initial data $(u^{0},u^{1})$ are in $\mathcal{D}(\mathcal{A})$, where we denote by
$$u^{0}=\left\{\begin{array}{ll}
u_{1}^{0}&\text{in }(0,\xi)
\\
u_{2}^{0}&\text{in }(\xi,1)
\end{array}\right.
\text{ and }
u^{1}=\left\{\begin{array}{ll}
u_{1}^{1}&\text{in }(0,\xi)
\\
u_{2}^{1}&\text{in }(\xi,1).
\end{array}\right.$$

The strong stability of energy for the model~\eqref{p1}-\eqref{p3} is provided if and only if $\xi$ is an irrational number (see~\cite{CCW} and~\cite{Ho}). Furthermore, for any $\xi\in(0,1)\backslash\Q$ the decay of the solution is not uniform in the energy space. For non symmetric boundary conditions (i.e Dirichlet boundary condition on one side and Newman boundary condition on the other side) the uniform exponential stability holds if and only if $\ds \xi=\frac{p}{q}$ with $p$ is odd (see~\cite{Ho}) and where the fastest decay rate of the solution is obtained when the actuator is located at the middle of the string (see~\cite{AHT}). Besides, if $\xi$ satisfies a Diophantine approximations properties then we have polynomial decay rate for the regular data (see~\cite{Tuc},~\cite{AHT} and~\cite{JTZ}). In our case of symmetric boundary conditions (Dirichlet boundary condition on both sides) Tucsnak~\cite{Tuc} proved that for every $\xi\in(0,1)\backslash\Q$ there exists $\ds\psi_{\xi}:[0,+\infty)\longrightarrow\R$ with $\ds\lim_{t\longrightarrow+\infty}\psi_{\xi}(t)=0$ such that the solution $(u_{1},u_{2})$ of~\eqref{p1}-\eqref{p3} satisfies the estimate
$$
\mathcal{E}(u_{1},u_{2})(t)\leq\psi_{\xi}(t).\|(u^{0},u^{1})\|_{\mathcal{D}(\mathcal{A})}^{2}\qquad\forall\,(u^{0},u^{1})\in\mathcal{D}(\mathcal{A}),\;t\geq 0
$$
where $\psi_{\xi}$ tends to zero at most as $\ds\frac{1}{t}$. At this stage we wondered if $\psi_{\xi}$ could tends to zero at least as $\ds\frac{1}{\ln(t)}$. The answer to this question is now given in the following main result, but first let's set
\begin{equation}\label{p19}
\mathcal{M}=\left\{\xi\in(0,1):\exists\,K_{1},\,K_{2}> 0,\;\left(\sin^{2}(\mu)+\sin^{2}(\xi\mu).\sin^{2}((1-\xi)\mu)\right)\e^{K_{1}\mu}\geq K_{2},\,\forall \mu\gg 1\right\},
\end{equation}
\begin{thm}\label{p4}
For any irrational number $\xi\in\mathcal{M}$ and for any $n\in\N$ there exists a constant $C_{\xi}>0$, such that for every initial data $(u^{0},u^{1})\in\mathcal{D}(\mathcal{A}^{n})$ the energy of the solution $(u_{1},u_{2})$ of~\eqref{p1}-\eqref{p3} satisfies
\[
\mathcal{E}(u_{1},u_{2})(t)\leq\frac{C_{\xi}}{(\ln(2+t))^{2n}}\|(u^{0},u^{1})\|_{\mathcal{D}(\mathcal{A}^{n})}^{2},\qquad\forall\,t>0.
\]
\end{thm}
It's clear that the assumption~\eqref{p19} include the set of $\xi\notin\Q$ satisfying the strong stability property which means that the energy is decreasing to zero as time goes to infinity. Moreover, we will show in section~\ref{p27} that this assumption is weaker than those given in~\cite{Tuc},~\cite{JTZ} and~\cite{AHT} where the polynomial stabilization is given depend on the Diophantine approximations properties of $\xi$. The above theorem is a consequence of Burq's result~\cite{Bur} which gives a sufficient condition to the resolvent estimate to obtain a decay rate of energy as given in Theorem~\ref{p4}. The main ingredient to prove the resolvent estimate is the use of what's called the technique of Carleman estimate. These kind of estimates have been used by several others to establish the logarithmic decay estimate for the dissipative systems (for instance in~\cite{LR},~\cite{B} and~\cite{I}), but only for multidimensional space problems. In our knowledge the Carleman estimates technique have never been used until now for one-dimensional space systems for the stabilization problems. However, in~\cite{BCV} this kind of estimates have been introduced for an inverse  problem. Besides, global Carleman estimates have been introduced for control problems in~\cite{BDL} and~\cite{Le} for parabolic equations.

This paper is organized as follows. In section~\ref{p5} we perform a suitable Carleman estimate. In section~\ref{p6} an appropriate resolvent estimate is established to prove Theorem~\ref{p4}. In Appendix~\ref{p28} some properties will be discussed on the assumption~\eqref{p19} to give it more meaning.
\section{Carleman estimate}\label{p5}
Let $[a,b]$ be an interval ($a<b$) in which we define the operator
\[
P=\frac{\ud^{2}}{\ud x^{2}}+\frac{1}{h^{2}},
\]
where $h>0$ is a small parameter. Let $\vfi\in\mathscr{C}^{4}([a,b])$, and define an adjoint operator $P_{\vfi}$ by
\[
P_{\vfi}=-h^{2}\e^{\vfi/h}P\e^{-\vfi/h}
\]
which can be written as follows
\[
P_{\vfi}w=(\de+i\vfi')^{2}w-w
\]
where we have denoted by $\ds\de=\frac{h}{i}\frac{\ud}{\ud x}$. The adjoint operator of $P_{\vfi}$ is given by
\[
P_{\vfi}^{*}w=(\de-i\vfi')^{2}w-w.
\]
We write $P_{\vfi}$ and $P_{\vfi}^{*}$ as the sum of a self-adjoint and an anti-adjoint operator as follows
\[
P_{\vfi}=Q_{2}+iQ_{1}\qquad\text{ and }\qquad P_{\vfi}^{*}=Q_{2}-iQ_{1}
\]
where $$Q_{2}=\frac{P_{\vfi}+P_{\vfi}^{*}}{2}=\de^{2}-(\vfi')^{2}-1\quad\text{ and }\quad Q_{1}=\frac{P_{\vfi}+P_{\vfi}^{*}}{2}=2\vfi'\de-ih\vfi''.$$
We perform some elementary calculations then by integration by parts we have
\begin{equation*}
\int_{a}^{b}v.Q_{2}\overline{w}\ud x=\int_{a}^{b}Q_{2}v.\overline{w}\ud x+ih\left(v(b).\de\overline{w}(b)+\de v(b).\overline{w}(b)-v(a).\de\overline{w}(a)-\de v(a).\overline{w}(a)\right)
\end{equation*}
and
\begin{equation*}
\int_{a}^{b}v.Q_{1}\overline{w}\ud x=\int_{a}^{b}Q_{1}v.\overline{w}\,\ud x+2ih\left(\vfi'(b).v(b).\overline{w}(b)-\vfi'(a).v(a).\overline{w}(a)\right).
\end{equation*}
This gives that
\begin{equation}\label{p7}
\int_{a}^{b}|P_{\vfi}w|^{2}\ud x=\int_{a}^{b}|Q_{2}w|^{2}\ud x+\int_{a}^{b}|Q_{1}w|^{2}\ud x+i\int_{a}^{b}[Q_{2},Q_{1}]w.\overline{w}\,\ud x+h\mathcal{B}(w),
\end{equation}
where the commutator bracket is defined by $[Q_{2},Q_{1}]=Q_{2}Q_{1}-Q_{1}Q_{2}$ and
\begin{equation*}
\begin{split}
\mathcal{B}(w)=2\vfi'(b)Q_{2}w(b).\overline{w}(b)-Q_{1}w(b).\de\overline{w}(b)-\de Q_{1}w(b).\overline{w}(b)
\\
-2\vfi'(a)Q_{2}w(a).\overline{w}(a)+Q_{1}w(a).\de\overline{w}(a)+\de Q_{1}w(a).\overline{w}(a).
\end{split}
\end{equation*}
The commutator $i[Q_{2},Q_{1}]$ can be written as follow
\[
i[Q_{2},Q_{1}]=h\left(4\vfi''\de^{2}-4ih\vfi'''\de+4\vfi''(\vfi')^{2}-h^{2}\vfi''''\right).
\]
Then a straightforward calculation leads to
\begin{equation}\label{p8}
\begin{split}
\re\left(i\int_{a}^{b}[Q_{2},Q_{1}]w.\overline{w}\,\ud x\right)=4h\int_{a}^{b}\vfi''|\de w|^{2}\ud x+4h\int_{a}^{b}\vfi''(\vfi')|w|^{2}\ud x
\\
-h^{3}\int_{a}^{b}\vfi''''|w|^{2}\ud x+4h^{2}\re\left(\frac{\vfi(b)}{i}\de w(b).\overline{w}(b)-\frac{\vfi(a)}{i}\de w(a).\overline{w}(a)\right).
\end{split}
\end{equation}
Similarly we also have
\begin{equation}\label{p9}
\begin{split}
\int_{a}^{b}|Q_{1}w|^{2}\ud x=4\int_{a}^{b}(\vfi')^{2}|\de w|^{2}\ud x+h^{2}\int_{a}^{b}(\vfi'')^{2}|w|^{2}\ud x+4h\re\left(\frac{1}{i}\int_{a}^{b}\vfi'\vfi''w.\de\overline{w}\,\ud x\right).
\end{split}
\end{equation}
We substitute~\eqref{p8} and~\eqref{p9} into~\eqref{p7} then we get
\begin{equation}\label{p10}
\begin{split}
\int_{a}^{b}|P_{\vfi}w|^{2}\ud x=\int_{a}^{b}|Q_{2}w|^{2}\ud x+4\int_{a}^{b}(\vfi')^{2}|\de w|^{2}\ud x+4h\int_{a}^{b}\vfi''|\de w|^{2}\ud x+h^{2}\int_{a}^{b}(\vfi'')^{2}|w|^{2}\ud x
\\
+4h\int_{a}^{b}\vfi''(\vfi')|w|^{2}\ud x-h^{3}\int_{a}^{b}\vfi''''|w|^{2}\ud x+4h\im\left(\int_{a}^{b}\vfi'\vfi''w.\de\overline{w}\,\ud x\right)
\\
-2h\vfi'(b)|\de w(b)|^{2}-2h^{2}\im(w(b).\de\overline{w}(b))-h(2\vfi'(b)(1+(\vfi'(b))^{2})-h^{2}\vfi'''(b))|w(b)|^{2}
\\
+2h\vfi'(a)|\de w(a)|^{2}+2h^{2}\im(w(a).\de\overline{w}(a))+h(2\vfi'(a)(1+(\vfi'(a))^{2})-h^{2}\vfi'''(a))|w(a)|^{2}.
\end{split}
\end{equation}
We assume that $\vfi$ satisfies the following three assumptions
\begin{itemize}
	\item[i)] $|\vfi'(x)|>0$ for every $x\in[a,b]$,
	\item[ii)] $\vfi''(x)>0$ for every $x\in[a,b]$,
	\item[iii)] $\vfi'(a)>0$.
\end{itemize}
The Carleman estimate is given by the following
\begin{thm}
Under the above assumptions on the weight function $\vfi$, there exist $h_{0}>0$ and $C>0$ such that for any $u\in H^{2}([a,b])$ with $u(a)=0$ we have
\begin{equation}\label{p11}
\begin{split}
h\int_{a}^{b}\e^{2\vfi/h}|u|^{2}\ud x+h^{3}\int_{a}^{b}\e^{2\vfi/h}|u'|^{2}\ud x+h^{3}|u'(a)|^{2}\e^{2\vfi(a)/h}
\\
\leq C\left(h^{4}\int_{a}^{b}\e^{2\vfi/h}|Pu|^{2}\ud x+\left(h|u(b)|^{2}+h^{3}|u'(b)|^{2}\right)\e^{2\vfi(b)/h}\right),
\end{split}
\end{equation}
for every $h\in(0,h_{0})$.
\end{thm}
\begin{pr}
Follow to~\eqref{p10} and weight function's assumption we have
\begin{equation*}
\begin{split}
\int_{a}^{b}|P_{\vfi}w|^{2}\ud x\geq C\bigg((1+h)\int_{a}^{b}|\de w|^{2}\ud x+h\int_{a}^{b}|w|^{2}\ud x-4h\left|\int_{a}^{b}\vfi'\vfi''w\de\overline{w}\ud x\right|
\\
+h|\de w(a)|^{2}+h|w(a)|^{2}-h^{2}|w(a).\de\overline{w}(a)|-h|\de w(b)|^{2}-h|w(b)|^{2}-h^{2}|w(b).\de\overline{w}(b)|\bigg).
\end{split}
\end{equation*}
By Young inequality and for $h$ small enough we follow that
\begin{equation}\label{p12}
\begin{split}
\int_{a}^{b}|P_{\vfi}w|^{2}\ud x\geq C\bigg(h\int_{a}^{b}|\de w|^{2}\ud x+h\int_{a}^{b}|w|^{2}\ud x+h|\de w(a)|^{2}+h|w(a)|^{2}
\\
-h|\de w(b)|^{2}-h|w(b)|^{2}\bigg).
\end{split}
\end{equation}
We set now $w=u\e^{\vfi/h}$ then we show easily that
$$
w'=u'\e^{\vfi/h}+h^{-1}\vfi'u\e^{\vfi/h}\text{ and }P_{\vfi}w=-h^{2}\e^{\vfi/h}Pu,
$$
and this leads to
$$
h\int_{a}^{b}|\de w|^{2}\ud x+h\int_{a}^{b}|w|^{2}\ud x\geq C\left(h^{3}\int_{a}^{b}\e^{2\vfi/h}|u'|^{2}\ud x+h\int_{a}^{b}\e^{2\vfi/h}|u|^{2}\ud x\right),
$$
and
$$
h|\de w(b)|^{2}+h|w(b)|^{2}\leq C(h^{3}|u'(b)|^{2}+h|w(b)|^{2})\e^{2\vfi(b)/h}.
$$
Hence by putting these inequalities into~\eqref{p12} and using the fact that $u(a)=0$ we can easily obtain~\eqref{p11}, and this achieves our proof.
\end{pr}

In terms of classical Carleman estimate and microlocal analysis, the so called subellipticity condition (see~\cite{LeRo2} and~\cite{LR}) given by
$$
\forall (x,\xi)\in[a,b]\times\R;\;p_{\vfi}(x,\xi)=0\,\Longrightarrow\,\{q_{2},q_{1}\}(x,\xi)\geq C>0,
$$
where the Poisson bracket is defined by $\{q_{2},q_{1}\}(x,\xi)=(\partial_{\xi}q_{2}\partial_{x}q_{1}-\partial_{x}q_{2}\partial_{\xi}q_{1})(x,\xi)$ and $p_{\vfi}=q_{2}+iq_{1}$, $q_{2}=\xi^{2}-((\vfi')^{2}+1)$ and $q_{1}=2(\vfi')\xi$ are respectively the principal symbol of the operators $P_{\vfi}$, $Q_{2}$ and $Q_{1}$, is verified in our case since $\p_{\vfi}(x,\xi)\neq0$ for all $(x,\xi)\in[a,b]\times\R$ thanks to the assumption $|\vfi'(x)|>0$.

If we have $u(b)=0$ instead of the Dirichlet condition $u(a)=0$, then by a simple change of variable $x=b+a-t$ the Carleman estimate holds just by permuting the role of the traces of $u$ in $a$ and $b$ in~\eqref{p11} and of course we have to replace the assumption $\vfi'(a)>0$ by $\vfi'(b)<0$.
\section{Resolvent estimate}\label{p6}
This section is devoted to establish the resolvent estimate, namely we prove that
\begin{equation}\label{p22}
\left\|(\mathcal{A}-i\mu\id)^{-1}\right\|\leq C\e^{K\mu},\qquad\mu\gg 1,
\end{equation}
which by~\cite{Bur} and~\cite{D} leads to the kind of decay energy estimate of the solution of system~\eqref{p1}-\eqref{p3} given in Theorem~\ref{p4}.

Let $(u_{1},u_{2},v_{1},v_{2})\in\mathcal{D}(\mathcal{A})$ and $(f_{1},f_{2},g_{1},g_{2})\in\mathcal{H}$ such that
$$
(\mathcal{A}-i\mu\id)\left(\begin{array}{c}
u_{1}
\\
u_{2}
\\
v_{1}
\\
v_{2}
\end{array}\right)=\left(\begin{array}{c}
f_{1}
\\
f_{2}
\\
g_{1}
\\
g_{2}
\end{array}\right),
$$
where $\mu$ is a real positive number large enough. 

This can be written in the following form
\begin{equation*}
\left\{\begin{array}{lll}
-i\mu u_{1}+v_{1}=f_{1}&\text{in}&(0,\xi)
\\
-i\mu u_{2}+v_{2}=f_{2}&\text{in}&(\xi,1)
\\
u_{1}''-i\mu v_{1}=g_{1}&\text{in}&(0,\xi)
\\
u_{2}''-i\mu v_{2}=g_{2}&\text{in}&(\xi,1)
\\
u_{1}(\xi)=u_{2}(\xi)&&
\\
u_{2}'(\xi)-u_{1}'(\xi)=v_{1}(\xi)&&
\\
u_{1}(0)=u_{2}(1)=0,&&
\end{array}\right.
\end{equation*}
which can be recast as the following boundary value problem
\begin{equation}\label{p13}
\left\{\begin{array}{lll}
v_{1}=f_{1}+i\mu u_{1}&\text{in}&(0,\xi)
\\
v_{2}=f_{2}+i\mu u_{2}&\text{in}&(\xi,1)
\\
u_{1}''+\mu^{2} u_{1}=g_{1}+i\mu f_{1}=\Phi_{1}&\text{in}&(0,\xi)
\\
u_{2}''+\mu^{2} u_{2}=g_{2}+i\mu f_{2}=\Phi_{2}&\text{in}&(\xi,1)
\\
u_{1}(\xi)=u_{2}(\xi)&&
\\
u_{2}'(\xi)-u_{1}'(\xi)=v_{1}(\xi)=f_{1}(\xi)+i\mu u_{1}(\xi)&&
\\
u_{1}(0)=u_{2}(1)=0.&&
\end{array}\right.
\end{equation}
Multiplying the third equation of~\eqref{p13} by $\overline{u}_{1}$ and the fourth one by $\overline{u}_{2}$ and integrating respectively over $(0,\xi)$ and $(\xi,1)$. Summing these two integrals then we show
\begin{equation}\label{p14}
\begin{split}
\int_{0}^{\xi}\Phi_{1}.\overline{u}_{1}\,\ud x+\int_{\xi}^{1}\Phi_{2}.\overline{u}_{2}\,\ud x=\mu^{2}\left(\int_{0}^{\xi}|u_{1}|^{2}\ud x+\int_{\xi}^{1}|u_{2}|^{2}\ud x\right)
\\
-\left(\int_{0}^{\xi}|u_{1}'|^{2}\ud x+\int_{\xi}^{1}|u_{2}'|^{2}\ud x\right)-i\mu|u_{1}(\xi)|^{2}-f_{1}(\xi).
\overline{u}_{1}(\xi).
\end{split}
\end{equation}
Taking the imaginary part of~\eqref{p14}, we obtain
\begin{equation*}
\begin{split}
\mu.|u_{1}(\xi)|^{2}\leq |f_{1}(\xi)|.|u_{1}(\xi)|+\left(\int_{0}^{\xi}|\Phi_{1}|^{2}\ud x\right)^{\frac{1}{2}}.\left(\int_{0}^{\xi}|u_{1}|^{2}\ud x\right)^{\frac{1}{2}}
\\
+\left(\int_{\xi}^{1}|\Phi_{2}|^{2}\ud x\right)^{\frac{1}{2}}.\left(\int_{\xi}^{1}|u_{2}|^{2}\ud x\right)^{\frac{1}{2}},
\end{split}
\end{equation*}
where by young inequality we follow that
\begin{equation}\label{p15}
\begin{split}
\mu.|u_{1}(\xi)|^{2}\leq C\Bigg(|f_{1}(\xi)|^{2}+\left(\int_{0}^{\xi}|\Phi_{1}|^{2}\ud x\right)^{\frac{1}{2}}.\left(\int_{0}^{\xi}|u_{1}|^{2}\ud x\right)^{\frac{1}{2}}
\\
+\left(\int_{\xi}^{1}|\Phi_{2}|^{2}\ud x\right)^{\frac{1}{2}}.\left(\int_{\xi}^{1}|u_{2}|^{2}\ud x\right)^{\frac{1}{2}}\Bigg).
\end{split}
\end{equation}

By solving the resolvent equation~\eqref{p13}, we find
\begin{equation*}
u_{1}(x)=\lambda_{1}\sin(\mu x)+\frac{1}{\mu}\int_{0}^{x}\sin(\mu(x-t))\Phi_{1}(t)\ud t
\end{equation*}
and
\begin{equation*}
u_{2}(x)=\lambda_{2}\sin(\mu(x-1))+\frac{1}{\mu}\int_{1}^{x}\sin(\mu(x-t))\Phi_{2}(t)\ud t,
\end{equation*}
where
\begin{equation*}
\begin{split}
\lambda_{1}=\frac{-\sin(\mu)+i\sin(\mu.\xi).\sin(\mu(1-\xi))}{\sin^{2}(\mu)+\sin^{2}(\mu.\xi).\sin^{2}(\mu(1-\xi))}\Bigg[\frac{\cos(\mu(1-\xi))}{\mu}\bigg(\int_{0}^{\xi}\sin(\mu(\xi-t))\Phi_{1}(t)\,\ud t
\\
+\int_{\xi}^{1}\sin(\mu(\xi-t))\Phi_{2}(t)\,\ud t\bigg)+\frac{\sin(\mu(1-\xi))}{\mu}\bigg(\int_{0}^{\xi}\e^{i\mu(\xi-t)}\Phi_{1}(t)\,\ud t
\\
+\int_{\xi}^{1}\cos(\mu(\xi-t))\Phi_{2}(t)\,\ud t+f_{1}(\xi)\bigg)\Bigg]
\end{split}
\end{equation*}
and
\begin{equation*}
\begin{split}
\lambda_{2}=\frac{-\sin(\mu)+i\sin(\mu.\xi).\sin(\mu(1-\xi))}{\sin^{2}(\mu)+\sin^{2}(\mu.\xi).\sin^{2}(\mu(1-\xi))}\Bigg[\frac{\cos(\mu\xi)+i\sin(\mu(1-\xi))}{\mu}
\\
\times\bigg(\int_{0}^{\xi}\sin(\mu(\xi-t))\Phi_{1}(t)\,\ud t+\int_{\xi}^{1}\sin(\mu(\xi-t))\Phi_{2}(t)\,\ud t\bigg)
\\
-\frac{\sin(\mu\xi)}{\mu}\bigg(\int_{0}^{\xi}\e^{i\mu(\xi-t)}\Phi_{1}(t)\,\ud t+\int_{\xi}^{1}\cos(\mu(\xi-t))\Phi_{2}(t)\,\ud t+f_{1}(\xi)\bigg)\Bigg].
\end{split}
\end{equation*}
Thus the derivative of $u_{1}$ and $u_{2}$ in $\xi$ are given by
\begin{equation}\label{p16}
u_{1}'(\xi)=\mu\lambda_{1}\cos(\mu.\xi)+\int_{0}^{\xi}\cos(\mu(\xi-t))\Phi_{1}(t)\ud t
\end{equation}
and
\begin{equation}\label{p17}
u_{2}'(\xi)=\mu\lambda_{2}\cos(\mu(\xi-1))+\int_{1}^{\xi}\cos(\mu(\xi-t))\Phi_{2}(t)\ud t.
\end{equation}

Let now $\vfi_{1}$ and $\vfi_{2}$ two weight functions defined respectively in $[0,\xi]$ and $[\xi,1]$ such that $\vfi_{1}\in\mathrm{C}^{4}([0,\xi])$ and $\vfi_{2}\in\mathrm{C}^{4}([\xi,1])$ and verifying
\begin{itemize}
	\item[i)] $|\vfi_{1}'(x)|>0$ in $[0,\xi]$ and $|\vfi_{2}'(x)|>0 $ in $[\xi,1]$
	\item[ii)] $\vfi_{1}''(x)>0$ in $[0,\xi]$ and $\vfi_{2}''(x)>0$ in $[\xi,1]$
	\item[iii)] $\vfi_{1}'(0)>0$ and $\vfi_{2}'(1)<0$.
\end{itemize}
We apply Carleman estimate (with taking $\ds h=\frac{1}{\mu}$), given by~\eqref{p11}, in each of the interval $(0,\xi)$ and $(\xi,1)$, respectively to the functions $u_{1}$ and $u_{2}$ solution of~\eqref{p13} with weight functions respectively $\vfi_{1}$ and $\vfi_{2}$. The sum of the two estimates leads to the following one
\begin{equation}\label{p18}
\begin{split}
h\int_{0}^{\xi}\e^{2\vfi_{1}/h}|u_{1}|^{2}\ud x+h\int_{\xi}^{1}\e^{2\vfi_{2}/h}|u_{2}|^{2}\ud x+h^{3}\int_{0}^{\xi}\e^{2\vfi_{1}/h}|u_{1}'|^{2}\ud x
\\
+h^{3}\int_{\xi}^{1}\e^{2\vfi_{2}/h}|u_{2}'|^{2}\ud x\leq C\bigg(h^{4}\int_{0}^{\xi}\e^{2\vfi_{1}/h}|\Phi_{1}|^{2}\ud x+h^{4}\int_{\xi}^{1}\e^{2\vfi_{2}/h}|\Phi_{2}|^{2}\ud x
\\
+h|u_{1}(\xi)|^{2}\left(\e^{2\vfi_{1}(\xi)/h}+\e^{2\vfi_{2}(\xi)/h}\right)+h^{3}|u_{1}'(\xi)|^{2}\e^{2\vfi_{1}(\xi)/h}+h^{3}|u_{2}'(\xi)|^{2}\e^{2\vfi_{2}(\xi)/h}\bigg).
\end{split}
\end{equation}
Substitute the explicit expression of $u_{1}'(\xi)$ and $u_{2}'(\xi)$ in~\eqref{p16} and~\eqref{p17} respectively into~\eqref{p18}, then by taking the maximum of $\vfi_{1}$ and $\vfi_{2}$ (that both of them are chosen strictly positive over the intervals $(0,\xi)$ and $(\xi,1)$ respectively) in the right hand side of~\eqref{p18} and their minimum in the left hand side, we find by the use of the assumption~\eqref{p19} for $h>0$ small enough that
\begin{equation}\label{p20}
\begin{split}
\int_{0}^{\xi}|u_{1}|^{2}\ud x+\int_{\xi}^{1}|u_{2}|^{2}\ud x+\int_{0}^{\xi}|u_{1}'|^{2}\ud x+\int_{\xi}^{1}|u_{2}'|^{2}\ud x
\\
\leq C\e^{C\mu}\bigg(\int_{0}^{\xi}|\Phi_{1}|^{2}\ud x+\int_{\xi}^{1}|\Phi_{2}|^{2}\ud x+|u_{1}(\xi)|^{2}\bigg).
\end{split}
\end{equation}
Combining~\eqref{p15} and~\eqref{p20}, one has
\begin{equation}\label{p21}
\begin{split}
\int_{0}^{\xi}|u_{1}|^{2}\ud x+\int_{\xi}^{1}|u_{2}|^{2}\ud x+\int_{0}^{\xi}|u_{1}'|^{2}\ud x+\int_{\xi}^{1}|u_{2}'|^{2}\ud x\leq C\e^{C\mu}\bigg(\int_{0}^{\xi}|\Phi_{1}|^{2}\ud x+\int_{\xi}^{1}|\Phi_{2}|^{2}\ud x
\\
+\left(\int_{0}^{\xi}|\Phi_{1}|^{2}\ud x\right)^{\frac{1}{2}}.\left(\int_{0}^{\xi}|u_{1}|^{2}\ud x\right)^{\frac{1}{2}}+\left(\int_{\xi}^{1}|\Phi_{2}|^{2}\ud x\right)^{\frac{1}{2}}.\left(\int_{\xi}^{1}|u_{2}|^{2}\ud x\right)^{\frac{1}{2}}+|f_{1}(\xi)|^{2}\bigg).
\end{split}
\end{equation}
Using the Young inequality and the fact that $H^{1}(0,\xi)\hookrightarrow\mathrm{C}(0,\xi)$ then we arrive at
\begin{equation*}
\begin{split}
\int_{0}^{\xi}|u_{1}|^{2}\ud x+\int_{\xi}^{1}|u_{2}|^{2}\ud x+\int_{0}^{\xi}|u_{1}'|^{2}\ud x+\int_{\xi}^{1}|u_{2}'|^{2}\ud x\leq C\e^{C\mu}\bigg(\int_{0}^{\xi}|\Phi_{1}|^{2}\ud x
\\
+\int_{\xi}^{1}|\Phi_{2}|^{2}\ud x+\int_{0}^{\xi}|f_{1}|^{2}\ud x+\int_{0}^{\xi}|f_{1}'|^{2}\ud x\bigg).
\end{split}
\end{equation*}
Now we have just to remember the expressions of $\Phi_{1}$ and $\Phi_{2}$ in~\eqref{p13} to show that
\begin{equation*}
\begin{split}
\int_{0}^{\xi}|u_{1}|^{2}\ud x+\int_{0}^{\xi}|u_{1}'|^{2}\ud x+\int_{\xi}^{1}|u_{2}|^{2}\ud x+\int_{\xi}^{1}|u_{2}'|^{2}\ud x\leq C\e^{C\mu}\bigg(\int_{0}^{\xi}|f_{1}|^{2}\ud x+\int_{0}^{\xi}|f_{1}'|^{2}\ud x
\\
+\int_{\xi}^{1}|f_{2}|^{2}\ud x+\int_{\xi}^{1}|f_{2}'|^{2}\ud x+\int_{0}^{\xi}|g_{1}|^{2}\ud x+\int_{\xi}^{1}|g_{2}|^{2}\ud x
\bigg),
\end{split}
\end{equation*}
which obviously leads to~\eqref{p22} and hence achieve the proof.
\appendix
\section{Appendix}\label{p28}
This section is devoted to some comments related to the assumption~\eqref{p19}.

In the case of non symmetric boundary conditions of type $u_{1}(0)=u_{2}'(1)=0$, it is well known that the strong stability holds if and only if $\ds\xi\neq\frac{p}{q}$ where $p$ is even and $q$ is odd. Moreover, if the assumption
$$
\left(\cos^{2}(\mu)+\cos^{2}(\xi\mu).\sin^{2}((1-\xi)\mu)\right)\e^{K_{1}\mu}\geq K_{2}
$$
takes place, for some $K_{1},\,K_{2}>0$, instead of that of~\eqref{p19}, we also have the same decay rate as given in Theorem~\ref{p4}. What is interesting in this condition is that it includes again the assumption of strong stabilization of the solutions as it is described above.

The objective in what follows is to prove that the condition~\eqref{p19} is attainable for a nonempty set of $\xi\in(0,1)$. Suppose that there exist a sequence of $\mu_{n}\longrightarrow+\infty$ such that
$$
\sin^{2}(\mu_{n})+\sin^{2}(\xi\mu_{n}).\sin^{2}((1-\xi)\mu_{n})\,\longrightarrow\,0.
$$
Hence, we can find two increasing functions $\psi_{1},\psi_{2}:\N\longrightarrow\N$ with $\ds\lim_{n\to+\infty}\psi_{1}(n)=\lim_{n\to+\infty}\psi_{2}(n)\!=+\infty$ and $\epsilon_{1},\epsilon_{2}:\N\longrightarrow\R$ with $\ds\lim_{n\to+\infty}\epsilon_{1}(n)=\lim_{n\to+\infty}\epsilon_{2}(n)=0$ such that
$$
\mu_{n}=\psi_{1}(n)\pi+\epsilon_{1}(n)\qquad\text{and}\qquad\xi\mu_{n}=\psi_{2}(n)\pi+\epsilon_{2}(n).
$$
We set $\nu_{n}=\psi_{1}(n)\pi$ and since $\nu_{n}=\mu_{n}-\epsilon_{1}(n)$, we obtain
$$
\nu_{n}\xi=\psi_{2}(n)\pi+\epsilon_{2}(n)-\xi\epsilon_{1}(n)\qquad\text{and}\qquad\nu_{n}(1-\xi)=(\psi_{1}(n)-\psi_{2}(n))\pi-\epsilon_{2}(n)+\xi\epsilon_{1}(n).
$$
Hence, we have
\begin{equation}\label{p23}
\sin^{2}(\nu_{n}\xi).\sin^{2}(\nu_{n}(1-\xi))\sim(\epsilon_{2}(n)-\xi\epsilon_{1}(n))^{2}((\epsilon_{1}(n)-\epsilon_{2}(n))-(1-\xi)\epsilon_{1}(n))^{2}.
\end{equation}
Now we suppose that~\eqref{p19} do not hold for every $K_{1},\,K_{2}>0$. If $\mu_{n}$ (as given above) and $c_{n}$ are the sequences which make that holds true namely,
$$\lim_{n\to+\infty}(\sin^{2}(\mu_{n})+\sin^{2}(\xi\mu_{n}).\sin^{2}((1-\xi)\mu_{n}))\e^{c_{n}\mu_{n}}=0\qquad\text{and}\qquad c_{n}\longrightarrow+\infty.$$
Since,
$$
\sin^{2}(\mu_{n})\e^{c_{n}\mu_{n}}\sim\epsilon_{1}(n)^{2}\e^{c_{n}\mu_{n}}
$$
and
$$
\sin^{2}(\xi\mu_{n})\sin^{2}((1-\xi)\mu_{n})\e^{c_{n}\mu_{n}}\sim\epsilon_{2}(n)^{2}(\epsilon_{1}(n)-\epsilon_{2}(n))^{2}\e^{c_{n}\mu_{n}}
$$
then from~\eqref{p23} the sequence $\nu_{n}=\psi_{1}(n)\pi$ contradicts also~\eqref{p19} for the same sequence $c_{n}$ i.e.,
\begin{equation}\label{p25}
\lim_{n\to+\infty}(\sin^{2}(\nu_{n})+\sin^{2}(\xi\nu_{n}).\sin^{2}((1-\xi)\nu_{n}))\e^{c_{n}\mu_{n}}=\lim_{n\to+\infty}\sin^{2}(\xi\nu_{n}).\sin^{2}((1-\xi)\nu_{n})\e^{c_{n}\nu_{n}}=0.
\end{equation}
We may write $|||\rho|||$ for the distance between $\rho\in\R$ and the nearest integer and let $\xi\in(0,1)$ be an irrational number that satisfies
\begin{equation}\label{p24}
\phi(m)|||m\xi|||\geq\kappa>0\qquad\forall\,m\gg 1\text{ and }m\in\N^{*},
\end{equation}
where $\phi$ is a positive increasing function of a real variable. Which imply obviously that
$$
\pi\phi(m)|||m\xi|||\geq\kappa\pi\qquad\forall\,m\gg 1\text{ and }m\in\N^{*}.
$$
With such a choice of $\xi$ it follows that $\ds0\leq\frac{\kappa\pi}{\phi\circ\psi_{1}(n)}\leq|||\psi_{1}(n)\xi|||\pi\leq\frac{\pi}{2}$ for every $n\in\N$, therefore, we show that
\begin{equation*}
\begin{split}
\sin^{2}(\xi\nu_{n}).\sin^{2}((1-\xi)\nu_{n})\e^{c_{n}\nu_{n}}=\sin^{4}(\psi_{1}(n)\xi\pi)\e^{c_{n}\psi_{1}(n)\pi}&=\sin^{4}(|||\psi_{1}(n)\xi|||\pi)\e^{c_{n}\psi_{1}(n)\pi}
\\
&\geq \sin^{4}\left(\frac{\kappa\pi}{\phi\circ\psi_{1}(n)}\right)\e^{c_{n}\psi_{1}(n)\pi}.
\end{split}
\end{equation*}
We suppose now that $\ds\lim_{x\to+\infty}\phi(x)=+\infty$ then we obtain
\begin{equation}\label{p26}
\sin^{4}\left(\frac{\kappa\pi}{\phi\circ\psi_{1}(n)}\right)\e^{c_{n}\psi_{1}(n)\pi}\sim\frac{\kappa^{4}\pi^{4}}{(\phi\circ\psi_{1}(n))^{4}}\e^{c_{n}\psi_{1}(n)\pi}.
\end{equation}
If $\phi:x\longrightarrow x$ is the identity function (in this case $\xi$ is called number of constant type also said to have bounded partial quotients~\cite{La}) then the right hand side of~\eqref{p26} goes to $+\infty$ as $n\longrightarrow+\infty$, which contradicts~\eqref{p25} and this means that $\xi$ satisfies~\eqref{p19}. Note in this case that a stronger result has been proved in~\cite{JTZ} and~\cite{Tuc}, namely it was shown that the energy of the solution decreases in a polynomial decay rate.
\\
More general, if $\xi$ satisfies~\eqref{p24} with $\phi(x)=O(\e^{cx})$ near $+\infty$ for some $c>0$ then by the same reasoning as earlier we are finding that $\xi$ does satisfy to the condition~\eqref{p19}. The following theorem gives more meaning to the choice of a such $\phi$.
\begin{thm}\cite[Theorem 32]{K}\label{p27}
Suppose that $f(x)$ is a positive continuous function of a positive variable $x$ and that $xf(x)$ is a non-increasing function. Then, the inequality
$$
|q\alpha-p|<f(q)
$$
has, for almost all $\alpha$, only a finite number of solutions in integers $p$ and $q$ (with $q>0$) if the integral
$$
\int_{c}^{+\infty}f(t)\,\ud t
$$
converges for some positive $c$.
\end{thm}
In particular, on the basis of Theorem~\ref{p27}, the inequality~\eqref{p24} is true under the following assumptions, namely
$$
\int_{c}^{+\infty}\frac{1}{\phi(t)}\,\ud t
$$
converges for some $c>0$ and $\ds x\longmapsto\frac{x}{\phi(x)}$ is a non-increasing function (for instance with $\phi(x)=x^{\alpha}\ln(x)^{(1+\epsilon)}$, for every constant $\epsilon>0$ and $\alpha\geq 1$ or also $\phi(x)=\e^{\beta x}$ for any $\beta>0$), for some $\kappa>0$ and for almost all $\xi\in(0,1)$.

For the polynomial decay we can find the same results as given in~\cite{JTZ}, in fact by using the resolvent method (see~\cite{BT}) where we proceed by the classical contradiction argument in which we explicit the solution of the resolvent problem as done in section~\ref{p6} then it follows under the assumption
\begin{equation}\label{p29}
\left(\sin^{2}(\mu)+\sin^{2}(\xi\mu).\sin^{2}((1-\xi)\mu)\right)\mu^{1+\epsilon}\geq K \quad\forall\, \mu\gg1,
\end{equation}
for some $\epsilon\geq 0$ and $K>0$, that the energy decay as follow
$$
\mathcal{E}(u_{1},u_{2})(t)\leq\frac{C_{\xi}}{(1+t)^{\frac{1}{1+\epsilon}}}\|(u^{0},u^{1})\|_{\mathcal{D}(\mathcal{A})}^{2},\qquad\forall\,t>0.
$$
for regular data $(u^{0},u^{1})\in\mathcal{D}(\mathcal{A})$. Noting here that, by proceeding as above we can show that assumption~\eqref{p29} generalize those given in~\cite{JTZ}, namely $\xi\in\mathscr{S}$  or $\xi\in B_{\epsilon}$.
\subsubsection*{Acknowledgments}
The author thanks the referees for many valuable remarks which helped us to improve the paper significantly. 
\nocite{*}
\bibliographystyle{alpha}
\bibliography{BibPWS}
\addcontentsline{toc}{section}{References}
\end{document}